\documentclass[12pt, oneside, a4paper]{article}%
\let\finishall\relax\let\Finishall\relax\let\getprepared\relax
\ifx\TestIngCommand\undefined\relax
\else\input ../../../proc/ppreamb-book.tex 
  \input init.tex \fi
\let\TestIngCommand\undefined

\usepackage{amsmath,amscd,amssymb}                                         
\usepackage{graphics}                                                      
\newtheorem{theo}{Theorem}                                                 
\newtheorem{lem}{Lemma}                                                    
\newtheorem{defi}{Definition}                                              
\newskip\ttglue\ttglue=.5em plus.25em minus.15em                           
\chardef\other=12                                                          
\def\ttverbatim{\begingroup                                                
  \catcode`\\=\other                                                       
  \catcode`\{=\other                                                       
  \catcode`\}=\other                                                       
  \catcode`\$=\other                                                       
  \catcode`\&=\other                                                       
  \catcode`\#=\other                                                       
  \catcode`\%=\other                                                       
  \catcode`\~=\other                                                       
  \catcode`\_=\other                                                       
  \catcode`\^=\other                                                       
  \obeyspaces \obeylines \tt}                                              
\outer\def\begintt{$$\let\par=\endgraf \ttverbatim \parskip=\z@            
  \catcode`\|=0 \rightskip-5pc \ttfinish}
{\catcode`\|=0 |catcode`|\=\other 
  |obeylines 
  |gdef|ttfinish#1^^M#2\endtt{#1|vbox{#2}|endgroup$$}}
\catcode`\|=\active                                                        
{\obeylines\gdef|{\ttverbatim\spaceskip\ttglue \let^^M=\  \let|=\endgroup}}
\def\firstname#1{\def\FIRSTNAME{#1}\ignorespaces}
\def\lastname#1{\def\LASTNAME{#1}\ignorespaces}
\def\middleinitial#1{\def\MIDDLEINI{#1}\ignorespaces}
\def\department#1{\def\DEPARTMENT{#1}\ignorespaces}
\def\institute#1{\def\INSTITUTE{#1}\ignorespaces}
\def\address#1{\def\ADDRESS{#1}\ignorespaces}
\def\country#1{\def\COUNTRY{#1}\ignorespaces}
\def\otheraffiliation#1{\def\OTHERAFFILIATION{#1}\ignorespaces}
\def\email#1{\def\EMAIL{#1}\ignorespaces}
\newcount\autcount\autcount=0                                              
\newcount\affcount\affcount=0                                              
\newcount\numcount\newcount\nummcount                                      
\newcount\nummmcount\newcount\nummmmcount                                  
\def\writename#1#2{\ \kern-1ex\hbox{
  \csname AUthor\the#1\endcsname\                                          
  \edef\TESTSTR{}\expandafter\ifx\csname auTHor\the#1\endcsname\TESTSTR    
  \else\csname auTHor\the#1\endcsname.\ \fi                                
  \csname authOR\the#1\endcsname$^{\csname AFF\the#1\endcsname}$
  \expandafter\ifx\csname corr\number#1\endcsname\relax                    
  \else\thanks{Corresponding author.}\ \fi                                 
  }\ifnum#1<#2, \else\ \kern-1ex\fi}
\def\writeemail#1{
  \nummcount=0\relax\nummmcount=0\relax                                    
  \loop\ifnum\nummcount<\autcount\advance\nummcount by1\relax              
    {\expandafter\ifnum\csname AFF\the\nummcount\endcsname=#1\relax        
    \global\advance\nummmcount by1\fi}\repeat                              
  \nummcount=0\relax\nummmmcount=0\relax                                   
  \loop\ifnum\nummcount<\autcount\advance\nummcount by1\relax              
    {\expandafter\ifnum\csname AFF\the\nummcount\endcsname=#1\relax        
    \global\advance\nummmmcount by1\relax\def\blank{}\expandafter          
    \ifx\csname EMAIL\the\nummcount\endcsname\blank(no e-mail)
    \else\csname EMAIL\the\nummcount\endcsname                             
    \fi                                                                    
    \ifnum\nummmmcount<\nummmcount; \fi\fi}\repeat}
\long\def\BeginAuthorList#1\EndAuthorList{#1\relax                         
  \author{\vbox{\hsize=390pt\noindent\numcount=0\relax                     
    \loop\ifnum\numcount<\autcount\advance\numcount by1\relax              
      \writename{\numcount}{\autcount}
      \repeat}\\[2mm]                                                      
    \vbox{\small\numcount=0\relax                                          
      \loop\ifnum\numcount<\affcount\advance\numcount by1\relax            
        \vbox{{\count0=\numcount\relax                                     
          \loop\expandafter\ifnum\csname AFF\the\count0\endcsname
            <\numcount\relax\advance\count0 by1\relax\repeat               
          $^{\csname AFF\the\count0\endcsname}$}
        \def\BLANK{}\expandafter\ifx\csname DEPT\the\numcount\endcsname    
          \BLANK                                                           
          \else\csname DEPT\the\numcount\endcsname, \fi                    
        \csname INST\the\numcount\endcsname,                               
        \csname ADDR\the\numcount\endcsname,                               
        \csname COUN\the\numcount\endcsname                                
        \edef\TEST{}\expandafter\ifx\csname OTHE\the\numcount\endcsname
          \TEST                                                            
          .\else;\break\csname OTHE\the\numcount\endcsname.\fi}
        \vbox{\writeemail{\numcount}}
        \repeat}\\}}
\expandafter\def\csname x1\endcsname{}
\expandafter\def\csname x2\endcsname{}
\expandafter\def\csname x3\endcsname{}
\expandafter\def\csname x4\endcsname{}
\expandafter\def\csname x5\endcsname{}
\expandafter\def\csname x6\endcsname{}
\expandafter\def\csname x7\endcsname{}
\expandafter\def\csname x8\endcsname{}
\expandafter\def\csname x9\endcsname{}
\def\Author#1#2{\global\advance\autcount by1\relax#2                       
  \expandafter\edef\csname AUthor\the\autcount\endcsname{\FIRSTNAME}
  \expandafter\edef\csname auTHor\the\autcount\endcsname{\MIDDLEINI}
  \expandafter\edef\csname authOR\the\autcount\endcsname{\LASTNAME}
  \expandafter\edef\csname EMAIL\the\autcount\endcsname{\EMAIL}
  \let\tempera\"\def\"{\string\"}\expandafter\ifx\csname x\DEPARTMENT
    \endcsname\relax                                                       
    \global\advance\affcount by1\relax\let\"\tempera                       
    \expandafter\edef\csname DEPT\the\affcount\endcsname{\DEPARTMENT}
    \expandafter\edef\csname INST\the\affcount\endcsname{\INSTITUTE}
    \expandafter\edef\csname ADDR\the\affcount\endcsname{\ADDRESS}
    \expandafter\edef\csname COUN\the\affcount\endcsname{\COUNTRY}
    \expandafter\edef\csname OTHE\the\affcount\endcsname{\OTHERAFFILIATION}
    \expandafter\edef\csname AFF\the\autcount\endcsname{\the\affcount}
  \else\expandafter\edef\csname AFF\the\autcount\endcsname{\DEPARTMENT}
  \fi\let\"\tempera\ignorespaces}
\def\CorrespondingAuthor#1#2{
  \expandafter\xdef\csname corr\number#1\endcsname{cor}
  \Author#1{#2}}
\def\PaperTitle#1{\title{\bf#1}}
\def\Category#1{\ignorespaces}
\def\keywords#1{{\noindent \emph{Keywords:}                                
  \def\BLANK{}\def\TEST{#1}\ifx\BLANK\TEST(n/a).\else#1\fi}}
\getprepared                                                               
\begin{document}                                                           
\PaperTitle{Minimal Prime  Ideals of  Ore Extensions over
 Commutative Dedekind Domains}%
\Category{(Pure) Mathematics}
\Category{Applied Mathematics}
\Category{Statistics}
\date{}

\BeginAuthorList 
  \Author1{
    \firstname{Amir Kamal}
    \lastname{Amir}
    \middleinitial{}   
    \department{Algebra  Research Division,
    Faculty of Mathematics and Natural Sciences}
    \institute{Institut Teknologi Bandung (ITB)}
    \address{Jl. Ganesha 10 Bandung 40132}
    \country{Indonesia}
     \otheraffiliation{Department of Mathematics, Hasanuddin
    University, Makassar 90245, Indonesia}
    \email{amirkamalamir@yahoo.com}}
  \Author2{
    \firstname{Pudji}
    \lastname{Astuti}
    \middleinitial{}    
    \department{Algebra  Research Division,
    Faculty of Mathematics and Natural Sciences}
    \institute{Institut Teknologi Bandung (ITB)}
    \address{Jl. Ganesha 10 Bandung 40132}
    \country{Indonesia}
    \otheraffiliation{}
    \email{pudji@math.itb.ac.id}}
  \Author3{
    \firstname{Intan}
    \lastname{Muchtadi-Alamsyah}
    \middleinitial{}    
    \department{2}
    \email{ntan@math.itb.ac.id}}
\EndAuthorList 
\maketitle 
\thispagestyle{empty} 
\begin{abstract} 

Let $R=D[x;\sigma, \delta]$ be an Ore extension over a commutative
Dedekind domain $D$, where $\sigma$ is an automorphism on $D$. In the
case $\delta =0$ Marubayashi et. al. already investigated the class of
minimal prime ideals in term of their contraction on the coefficient ring $D$.
In this note we extend this result to a general case $\delta \ne 0$.

\end{abstract} 


\finishall 
\section{Introduction}

This paper studies minimal prime ideals of an Ore extension over
a commutative Dedekind domain. Ore extensions are widely used as
the underlying rings of various linear systems investigated in the area Algebraic
system theory. These systems may represent systems coming  from mathematical physics, applied mathematics
and engineering sciences which can be described by means of systems of
ordinary or partial differential equations, difference equations,
differential time-delay equations, etc. If these systems are linear,
they can be defined by means of matrices with entries in
non-commutative algebras of functional operators such as the ring of
differential operators, shift operators, time-delay operators, etc.
An important class of such algebras is called Ore extensions (Ore
Algebras).


The structure of prime ideals of various kind of Ore extensions have
been investigated during the last few years. In \cite{IRV2},
\cite{LMA}  primes of Ore extensions over commutative noetherian
rings were considered. In \cite{CHI}, \cite{MFM}, and \cite{DSP},
prime ideals of Ore extensions of derivation type were described.
These result recently were exploited in \cite{MAW} to investigate
properties of minimal prime rings of Ore extensions of derivation
type
 in term of their contraction on the coefficient ring.
In this note we extend this result to a general Ore extension of
automorphism type .


\section{Ore Extension}
We recall some definitions, notations, and more or less well known
facts concerning. A  $\textit{(left) skew derivation}$ on a ring $D$
is a pair $(\sigma, \delta)$ where $\sigma$ is a ring endomorphism
of $D$ and $\delta$ is a $\it{(left)}$
~$\sigma$-$\textit{derivation}$ ~on $D$; that is, an additive map
from $D$ to itself such that $\delta(ab) = \sigma(a) \delta(b) +
\delta(a)b$ for all $a,b \in D$. For $(\sigma, \delta)$ any  skew
derivation on a ring $D$, we obtain
\[
\delta(a^m) = \sum^{m-1}_{i=0} \sigma(a)^i \delta(a) a^{m-1-i}
\]
for all $a \in D$ and $m=1,2, \cdots $. (See \cite[Lemma 1.1]{GOO})

\begin{defi}
 Let $D$ be a ring with identity 1 and $(\sigma, \delta)$ be a
(left) skew derivation on the ring $D$.
 \noindent \textit{The Ore Extension}   over $D$ with respect to the skew derivation
$(\sigma, \delta)$
 is the ring consisting of all polynomials over $D$ with
 an indeterminate $x$ denoted by:
\[ D[x;\sigma ,\delta ]\, \, =\, \, \{ \, f(x)\, =\, a_{n}
x^{n} \, +\, \, \cdots \, \, +\, a_{0} \, \, \mid \, \, a_{i} \, \in
\, D\, \, \} \]
 satisfying the following equation, for all $a \in D$
\[xa\, \, =\, \sigma (a)x\, +\, \delta (a).\]
\end{defi}
The notations  $D[x; \sigma]$ and  $D[x; \delta]$ stand for the particular Ore extensions
where respectively $\delta = 0$ dan $\sigma$ the identity map.
For the case $\delta = 0$, Marubayashi et. al. \cite{MAW} studied the factor rings of
$D[x; \sigma]$ over minimal prime ideals where $D$ is a commutative Dedekind domain.
In order to extend their results to general cases, this paper investigates the class of
minimal prime ideals in $D[x; \sigma, \delta]$.

The Ore extension $D[x; \sigma, \delta]$
is a free left $D$-module with basis $1, x, x^2, \cdots$
To abbreviate the assertion, the symbol $R$ stands for
the Ore extension $D[x; \sigma, \delta]$ constructed from
a ring $D$ and a skew derivation $(\sigma, \delta)$ on $D$.
The $\textit{degree}$ of a nonzero element $f \in R$ is defined in
the obvious fashion. Since the standard form for elements of $R$ is
with left-hand coefficients, the $\textit{leading coefficient}$ of
$f$ is $f_n$ if
\[
f= f_0 + f_1x + \cdots + f_{n-1}x^{n-1} + f_nx^n
\]
with all $f_i \in D$ and $f_n \neq 0$. If $\sigma$ is an
automorphism, $f$ can also be written with right-hand coefficients,
but then its $x^n$-coefficient is $\sigma^{-n}(f_n)$.
While a general formula for $x^na$ where $a \in D$ and $n \in
\mathbb{N}$ is too involved to be of much use, an easy induction
establishes that
\[
x^na = \delta^n(a) + a_x + \cdots + a_{n-1}x^{n-1} + \sigma^n(a)x^n
\]
for some $a_1, \cdots , a_{n-1} \in D$.\\

In preparation for our analysis of the types of ideals occured
when  prime ideals of an Ore extension $D[x; \sigma, \delta]$ are
contracted to the coefficient ring $D$, we consider $\sigma$-prime,
$\delta$-prime, and $(\sigma, \delta)$-prime ideals of $D$.

\begin{defi}
Let $\Sigma$ be a set of map from the ring $D$ to itself. A
$\Sigma$-ideal of $D$ is any ideal $I$ of $D$ such that $\alpha(I)
\subseteq I$ for all $\alpha \in \Sigma$. A $\Sigma$-prime ideal is
any proper $\Sigma$-ideal $I$ such that whenever $J, K$ are
$\Sigma$-ideals satisfying $JK \subseteq I$, then either $J
\subseteq I$ or $K \subseteq I$.
\end{defi}
In the context of a ring $D$ equipped with a skew derivation
$(\sigma, \delta)$, we shall make use of the above definition in
the cases $\Sigma = \{\sigma \}, \Sigma = \{\delta \}$ or $\Sigma =
\{\sigma, \delta \}$; and simplify the prefix $\Sigma$ to respectively $\sigma,
\delta,$ or $(\sigma, \delta)$. Concerning the contraction of
prime ideals in an Ore Extension to its coefficient ring,
Goodearl \cite{GOO} obtained the following theorem which will be of use later.

\begin{theo}
Let $R=D[x; \sigma, \delta]$ where $D$ is a commutative Dedekind
domain and $\sigma$ is an automorphism. If $\mathfrak{p}$ is any
ideal of $D$ which is $(\sigma, \delta)$-prime, then $\mathfrak{p}=P
\cap R$ for some prime ideal $P$ of $R$ and more specially
$\mathfrak{p}R \in$ Spec($R$) where $Spec(R)$ denotes the set of all Prime ideal in $R$.
\end{theo}

\section{Minimal Prime Ideals of Ore Extensions}

Throughout this section, let $D$ be a commutative Dedekind domain and
$R = D[x;\sigma, \delta]$ be the Ore
extension over $D$, for  $(\sigma, \delta)$ is a skew derivation, $\sigma \neq 1$ is an automorphism of $D$ and $\delta \neq 0$.

Marubayashi et. al. \cite{MAW} already investigated the class of minimal prime ideals in $R$ where
$\delta=0$. For the case $\delta \ne 0$ but it is inner, Goodearl \cite{GOO} showed the
existance of  an isomorphism between $D[x;\sigma, \delta]$ and $D[y;\sigma]$ as described in the following.
The $\sigma$-derivative $\delta$ is called  inner if  there exists an element $a \in D$ such that $\delta(b) = ab
- \sigma(b)a$ for all $b \in D$.

\begin{theo} \label{t-isomorfisma}
Let $D[x;\sigma, \delta]$ be an Ore extension where $\sigma \neq 1$
and $\delta \neq 1$. If $\delta$ is an inner $\sigma$-derivation,
i.e, there exists $a \in D$ such that $\delta(b) = ab - \sigma(b)a$
for all $b \in D$,  then $D[x;\sigma,\delta]$ and  $D[y;\sigma ]$, where $y=x-a$, are isomorph.
\end{theo}

Hence, by combining Theorem \ref{t-isomorfisma} and the class of minimal prime ideals in Ore extensions $R$ obtained in \cite{MAW}
we can derive the class of minimal prime ideals in Ore extension $R=D[x; \sigma, \delta]$ for $\sigma$ and $\delta$ being
respectively an automorphism and inner as the following.
Notation $\textnormal{Spec}_0 (R)$ stands for the set of all prime ideals in $R= D[x;\sigma,\delta]$ having
zero intersection with $D$.

\begin{theo} \label{t-utkdeltainner}
Let $R=D[x;\sigma,\delta]$ be an Ore extension.\\
\textnormal{(1)} The set \\
$ \{ \mathfrak{p}[x;\sigma,\delta], ~P \mid \mathfrak{p}\ \hbox{is
a}\ \sigma-\hbox{prime ideal of}\ D \hbox{and}\ P \in
\textnormal{Spec}_0 (R) ~\hbox{with} ~P \neq (0) \}$
\\
consists  of all
minimal prime ideals of $R$.\\
\textnormal{(2)} Let $P \in \textnormal{Spec}(R)$ with $P \neq (0)$.
Then $P$ is invertible if and only if it is a minimal prime ideal of
$R$.
\end{theo}

Now we shall investigate the class of minimal prime ideals for general $\delta \ne 0$.
For this general case, we need the following lemma.

\begin{lem} \label{lema1}
If $P = \mathfrak{p}[x;\sigma,\delta]$ is a minimal prime ideal of
$R$ where $\mathfrak{p}$ is a  $(\sigma,\delta)$-prime ideal of
$D$, then $\mathfrak{p}$ is a minimal $(\sigma,\delta)$-prime ideal
of $D$.
\end{lem}
\textit{Proof}. Assume that $\mathfrak{p}$ is a
$(\sigma,\delta)$-prime ideal of $D$ but  is not  minimal
$(\sigma,\delta)$-prime. Let $\mathfrak{q}$ be a
$(\sigma,\delta)$-prime ideal of $D$ such that  $\mathfrak{q}
\varsubsetneq \mathfrak{p}$ and $\mathfrak{q}\ne (0)$. Then applying \cite[Theorem 3.1]{CHI} we have
$\mathfrak{q}R \in$ Spec($R$). So,
$\mathfrak{q}R \varsubsetneq \mathfrak{p}R = P$. This is a
contradiction because $P$ is a minimal prime ideal. $\square $

\begin{theo}
 Let $P$ be a prime ideal of $R$ and $ P \cap D = \mathfrak{p} \ne
(0)$. Then $P$ is a minimal prime ideal of $R$ if and only if either
$P =\mathfrak{p}[x;\sigma,\delta]$ where $\mathfrak{p}$ is a minimal
$(\sigma,\delta)$-prime ideal of $D$ or $(0)$ is the largest
$(\sigma,\delta)$-ideal of $D$ in $\mathfrak{p} $.
\end{theo}
\noindent \textit{Proof}.~Let $P$ be a minimal prime ideal of $R$
and $ P \cap D = \mathfrak{p} \ne (0)$.
Since $\mathfrak{p} \neq (0)$, then there are two cases \cite[Theorem 3.1]{GOO}; namely, either $\mathfrak{p}$ is a $(\sigma,
\delta)$-prime ideal of $D$ or $\mathfrak{p}$ is a prime ideal of
$D$ and
$\sigma(\mathfrak{p}) \ne \mathfrak{p}$.

If $\mathfrak{p}$ is a $(\sigma, \delta)$-prime ideal of $D$, then
$\mathfrak{p}R \in$ Spec($R$), by \cite[Theorem 3.1]{GOO}. So,
$\mathfrak{p}R = P$ because $\mathfrak{p}R \subseteq P$ and $P$ is a
minimal prime ideal. On the other hand $\mathfrak{p}R =
\mathfrak{p}[x; \sigma, \delta]$. From here, we get $P =
\mathfrak{p}[x; \sigma, \delta]$ and $\mathfrak{p}$ is a minimal
$(\sigma, \delta)$-prime ideal of $D$,
by Lemma \ref{lema1}.

Suppose $\mathfrak{p}$ is a prime ideal of $D$ and $\sigma(\mathfrak{p})
\ne \mathfrak{p}$. Let $\mathfrak{m}$ be the largest $(\sigma,
\delta)$-ideal contained in  $\mathfrak{p}$ and assume that
$\mathfrak{m} \ne (0)$. Then by primeness of $\mathfrak{p}$ it can be shown
 that $\mathfrak{m}$ is a $(\sigma, \delta)$-prime ideal of
$D$. So, $\mathfrak{m}R$ is a prime ideal of $R$ by \cite[Proposition 3.3]{GOO}. On the other hand, since $\sigma(\mathfrak{p}) \ne
\mathfrak{p}$, we have $\mathfrak{m} \varsubsetneq \mathfrak{p}$.
So, $\mathfrak{m}R \varsubsetneq \mathfrak{p}R \subseteq P$, i.e,
$P$ is not a minimal prime. This is a contradiction. So, $(0)$ is
the largest $(\sigma,\delta)$-ideal of $D$ in  $\mathfrak{p}. $

Conversely, let $ P = \mathfrak{p}[x; \sigma, \delta]$, where $\mathfrak{p}$
is a minimal $(\sigma, \delta)$-prime ideal of $D$. Then, according to  \cite[Theorem 3.3]{GOO}, $ P =
\mathfrak{p}[x; \sigma, \delta]$ is a prime ideal of $R$. Let $Q$ be a prime ideal of $R$ where $Q \subseteq P$.
Set $\mathfrak{q} = Q \cap D$, then $\mathfrak{q} = Q \cap D
\subseteq  P \cap D = \mathfrak{p}$. Applying \cite[Theorem 3.1]{GOO}, we have
two cases; namely, either $\mathfrak{q}$ is a $(\sigma,
\delta)$-prime ideal of $D$ or $\mathfrak{q}$ is a prime ideal of
$D$ and $\sigma(\mathfrak{q}) \ne \mathfrak{q}$.

For the first case, suppose $\mathfrak{q}$ is a $(\sigma, \delta)$-prime ideal of $D$. Then
$\mathfrak{q}=\mathfrak{p}$ because $\mathfrak{q}\subseteq
\mathfrak{p}$ and $\mathfrak{p}$ is a minimal $(\sigma,
\delta)$-prime ideal of $D$. So, $ P = \mathfrak{p}[x; \sigma,
\delta]= \mathfrak{q}[x; \sigma, \delta] \subseteq Q$. This implies
$P=Q$.
For the other case, if $\mathfrak{q}$ is a prime ideal of $D$, then
$\mathfrak{q}=\mathfrak{p}$ because $D$ is a Dedekind domain. So, $
P = \mathfrak{p}[x; \sigma, \delta]= \mathfrak{q}[x; \sigma, \delta]
\subseteq Q$. This implies $P=Q$. Thus
 $P$ is a minimal prime ideal of $R$.

Let $(0)$ be the largest $(\sigma, \delta)$-ideal of $D$ in
$\mathfrak{p}$. Let $Q$ be a prime nonzero  ideal of $R$ satisfying $Q \subseteq
P$. Set $\mathfrak{q} = Q \cap D$, then $\mathfrak{q} = Q \cap D
\subseteq  P \cap D = \mathfrak{p}$. Similar to the above explanation, we have
two cases; namely, either $\mathfrak{q}$ is a $(\sigma,
\delta)$-prime ideal of $D$ or $\mathfrak{q}$ is a prime ideal of
$D$ and $\sigma(\mathfrak{q}) \ne \mathfrak{q}$ but the first case will not happen. If it is so,
 then, because of $(0)$ being the largest $(\sigma,
\delta)$-ideal of $D$ in $\mathfrak{p}$,
$\mathfrak{q}= (0)$ implying a contradiction  $Q \cap D = 0$ (see  \cite[Lemma 2.19]{GOW}).
Thus $\mathfrak{q}$ is a prime ideal of $D$ with $\sigma(\mathfrak{q})
\ne \mathfrak{q}$. As a result, $\mathfrak{q}= \mathfrak{p}$. So $Q \cap D =
P \cap D$, which, according to   \cite[Proposition 3.5]{GOW}, implies $Q = P$.
Thus  $P$ is the
minimal prime ideal of $R$.  $\square $

\section{Concluding Remark}
In this paper we identify all minimal prime ideals in an Ore extension over
a Dedekind domain according to their contraction on the coefficient ring.
Studying the results in Marubayashi et. al \cite{MAW} it is  expected that this identification can be used to study the structure of the corresponding factor rings which is currently under investigation.

\section*{Acknowledgement}
This research was supported by ITB according to Surat Perjanjian Pelaksanaan Riset Nomor: 265/K01.7/PL/2009, 6 February 2009.

\Finishall 
\end{document}